\documentclass[number,citesort,dvips]{arxbj}
\usepackage{cursive} %uplgreek}

\aid{0}
\volume{16}
\issue{4}
\pubyear{2010}
\firstpage{1114}
\lastpage{1136}
\doi{10.3150/09-BEJ246}

\makeatletter
\newtheorem{theorem}{Theorem}[section]
\newtheorem{corollary}{Corollary}[section]
\newtheorem{lemma}{Lemma}[section]

\newremark{Remarks}{Remark}[section]

\makeatother

\begin{document}
\begin{frontmatter}

\title{Compound Poisson and signed compound Poisson approximations to the
Markov binomial law}
\runtitle{Compound Poisson--Markov binomial}

\begin{aug}
\author[a]{\fnms{V.} \snm{\v{C}ekanavi\v{c}ius}\thanksref{a}\ead[label=e1]{vydas.cekanavicius@mif.vu.lt}\corref{}} \and
\author[b]{\fnms{P.} \snm{Vellaisamy}\thanksref{b}\ead[label=e2]{pv@math.iitb.ac.in}}
\runauthor{V. \v{C}ekanavi\v{c}ius and P. Vellaisamy}
\address[a]{Department of Mathematics and Informatics, Vilnius University,
Naugarduko 24, Vilnius 03225, Lithuania. \printead{e1}}
\address[b]{Department of Mathematics, Indian Institute of Technology Bombay,
Powai, Mumbai 400076, India. \printead{e2}}
\end{aug}

% HISTORY:
\received{\smonth{7} \syear{2007}}
\revised{\smonth{10} \syear{2009}}

% ABSTRACT
%
\begin{abstract}
Compound Poisson distributions and signed compound Poisson measures
are used for approximation of the Markov binomial distribution. The
upper and lower bound estimates are obtained for the total
variation, local and Wasserstein norms. In a special case,
asymptotically sharp constants are calculated. For the upper bounds,
the smoothing properties of compound Poisson distributions are
applied. For the lower bound estimates, the characteristic function
method is used.
\end{abstract}

% KEYWORDS
%

\end{frontmatter}

%s1 ###
\section{Introduction}\label{sec1}
% Initial capital letter, then lower case. No
%full stop.

The closeness of a compound Poisson (CP) distribution to the Markov
binomial (MB) distribution has been investigated in numerous papers;
see, for
example, \cite{CM99,CM01,Dob53,Gan82,Koo50,Ser86,Ve99,Wan92} and the
references therein. Related problems were considered in
\cite{CGS94,MBB,CV02,Erh99,Aki93,Hsiau,Ser75,Vel04} and \cite{WCC03}.
One would expect the MB--CP case to have been comprehensively studied.
As it
turns out, this is not the case. Many papers deal with the convergence facts
only. Only a few of the papers dealing with the estimates of
accuracy of approximation involve no assumptions about the stationarity
of the Markov chain.

%With a very few exceptions the 'magic factor' is absent from the
%estimates of accuracy of approximation.
%
% As
%observed in [2, p. 5], in total variation Poisson approximation to
%the sum of independent indicator variables
% is sharper than the sum of local estimates.
%Indeed, approximation of Bernoulli variable by the Poisson one
%
% Some additional factor appears just like by 'magic'.
% One can expect similar factor to appear for the sum of dependent
% indicators. However, in MB--CP case,

The aim of this paper is to discuss some compound approximations
for non-stationary Markov chains. We show that for our version of
the MB distribution, the natural approximation is a convolution of
CP and compound binomial distributions, both having the same
compounding geometric law. We outline some principles of
construction of asymptotic expansions and consider second order
approximations. Part of the paper is devoted to signed compound
Poisson approximations which can be viewed as the second order
expansions {in the exponent}. We obtain upper and lower bound
estimates and show that under certain conditions, they are of the
same order of accuracy. All estimates are proved for the total
variation, local and Wasserstein norms.
% and contain 'magic' factors.
For the upper bound estimates, we employ a convolution technique which
can be dated back to \cite{Lec60}. For the lower bound estimates,
we use the characteristic function method. The methods of proof do
not allow for reasonably small absolute constants. However, in
special cases, asymptotically sharp constants are calculated.

We now introduce some notation. Let $I_k$ denote
the distribution concentrated at an integer $k\in\mathbb{Z},$ the set of
integers, and set $I=I_0$. In what follows, $V$ and
$M$ denote two finite signed measures on $\mathbb{Z}$. Products and
powers of
$V$ and $M$ are understood in the convolution sense, that is,
$VM\{A\} = \sum_{k=-\infty}^\infty V\{A-k\} M\{k\}$ for a set
$A\subseteq\mathbb{Z}$; further, $M^0=I$.
The total variation norm, the local norm and the
Wasserstein norm of~$M$ are denoted by
\[
\|M\|=\sum_{k=-\infty}^\infty|M\{k\}|,\qquad
\|M\|_\infty= \sup_{k\in\mathbb{Z}}|M\{k\}|,\qquad
\|M\|_{\mathrm{W}}=\sum_{k=-\infty}^\infty|M\{(-\infty,k]\}|,
\]
respectively. Note that $\|(I_1-I)M\|_{\mathrm{W}}=\|M\|$.
The logarithm and exponential of~$M$ are given, respectively, by
\[
\ln M=\sum_{k=1}^\infty\frac{(-1)^{k+1}}{k} (M-I)^k\qquad
(\mbox{if }\|M-I\|<1),\qquad
\mathrm{e}^M=\exp\{M\}=\sum_{k=0}^\infty\frac{1}{k!} M^k.
\]
Note that
\[
\|VM\|_\infty\leq\|V\|\|M\|_\infty,\qquad
\|VM\|\leq\|V\|\|M\|,\qquad
\|\mathrm{e}^{M}\|\leq\mathrm{e}^{\|M\|}. %\label{star}
\]
Let $\widehat{M}(t)$
$(t\in\mathbb{R})$ be the Fourier transform of $M$. We denote by $C$
positive absolute constants. $\Theta$ stands for any
finite signed measure on $\mathbb{Z}$ satisfying $\|\Theta\|\leq
1$. The values of $C$ and $\Theta$ can vary from line to line, or
even within the same line. Sometimes, to avoid possible ambiguity,
the $C$'s are supplied with indices. For $x\in\mathbb{R}$ and
$k\in\mathbb{N}=\{1,2,3,\ldots\}$, we set
\[
\pmatrix{x \cr k}=\frac{1}{k!} x(x-1)\cdots(x-k+1),\qquad
\pmatrix{x \cr0}=1.
\]
Let $\xi_0,\xi_1,\ldots,\xi_n,\ldots$ be a Markov
chain with the initial distribution
\[
\mathrm{P}(\xi_0=1)=p_0,\qquad \mathrm{P}(\xi_0=0)=1-p_0, \qquad
p_0\in[0,1]
\]
and transition probabilities
\begin{eqnarray*}
& \mathrm{P}(\xi_i=1 | \xi_{i-1}=1)=p,\qquad
\mathrm{P}(\xi_i=0 | \xi_{i-1}=1)=q,&\\
& \mathrm{P}(\xi_i=1 | \xi_{i-1}=0)=\overline{q},\qquad
\mathrm{P}(\xi_i=0 | \xi_{i-1}=0)=\overline{p},&\\
& p+q=\overline{q}+\overline{p}=1,\qquad
p,\overline{q}\in(0,1), \qquad
i\in\mathbb{N}.&
\end{eqnarray*}
The distribution of $S_n=\xi_1+\cdots+\xi_n$ $(n\in\mathbb{N})$ is called
the \textit{Markov binomial distribution}. We denote it by $F_n$. We should
note that the definition of the Markov binomial distribution
varies slightly from paper to paper; see \cite{Dob53,Ser75} and \cite
{Wan81}. We choose the definition which, on
the one hand, contains the binomial distribution as a special case
and, on the other hand, allows comparison to the Dobrushin's
results. Dobrushin \cite{Dob53} assumed that $p_0=1$ and considered
$S_{n-1}+1$.\looseness=1

Later, we will need various characteristics of $S_n$. Let
\begin{eqnarray*}
\gamma_1&=&\frac{q\overline{q}}{q+\overline{q}},\qquad
\gamma_2=-\frac{q\overline{q}^2}{(q+\overline{q})^2} \biggl
(p+\frac{q}{q+\overline{q}} \biggr)-\frac{\gamma_1^2}{2},
\\
\gamma_3&=&\gamma_1^2\tilde\gamma_3,
\\
\tilde\gamma_3&=&\frac{\gamma_1}{3}+ \frac{1}{q(q+\overline{q})}
\biggl\{
p^2\overline{q}+\frac{pq(2\overline{q}-q)}{q+\overline{q}}+\frac
{2\overline{q}
q^2}{(q+\overline{q})^2} \biggr\}+\frac{\overline{q}}{q+\overline
{q}} \biggl(p+\frac{q}{q+\overline{q}} \biggr),
\\
\lambda&=&n-p_0,\qquad\varkappa_1=\gamma_1 \biggl(\frac{\overline
{q}-p}{q+\overline{q}}-p_0 \biggr),\qquad
\varkappa_2=p_0\frac{pq}{q+\overline{q}},\qquad C_1=\ln\frac{30}{19}
=0.4567\ldots.
\end{eqnarray*}
We use the following measures also:
\begin{eqnarray*}
G&=&qI_1\sum_{j=0}^\infty p^jI_j\qquad
\biggl(\widehat G(t)=\frac{q\mathrm{e}^{\mathrm{i}t}}{1-p\mathrm
{e}^{\mathrm{i}t}} \biggr),\qquad
H=I+\varkappa_2(G-I),\\
H_1&=&(1-\gamma_1)I+\gamma_1 G,\qquad
H_1^\lambda=\exp\Biggl\{\lambda\sum_{j=1}^\infty\frac{(-1)^j}{j}\gamma
_1^j(G-I)^j \Biggr\}\qquad
\bigl(\widehat H_1^\lambda(t)=(\widehat H_1(t))^\lambda\bigr),\\
%D_1&=&\exponent{\gamma_1(G-\dirac)}, D_1^\lambda=\exponent{\lambda
% D_2&=&\exponent{\gamma_1(G-\dirac)+\gamma_2(G-\dirac)^2},
%D_3=\exponent{\gamma_1(G-\dirac)+\gamma_2(G-\dirac)^2+\gamma_3(G-
D_j&=&\exp\Biggl\{\sum_{i=1}^{j} \gamma_i(G-I)^{i}\Biggr\},\qquad
1\leq j \leq3, \qquad D_1^\lambda=\exp\{\lambda\gamma_1(G-I)\}.
\end{eqnarray*}

%s2 ###
\section{Known results}

In this section, we discuss some of the known results on the
compound Poisson approximations to the MB distribution. Many papers
deal with the convergence facts only; see, for example,
\cite{Gan82,Hsiau,Koo50,Wan81}.
%In the
%absolute majority of the papers dealing with the accuracy of CP
%approximations the estimates contain no `magic' factor.
Usually, the chain is assumed to be stationary. A typical example
is Theorem 4.1 in \cite{Ser86} which states that if
$p_0=\overline{q}/(q+\overline{q})$ and $\tilde\alpha>0$, then
%e1 ###
%
\begin{equation}
\|F_n-\exp\{\tilde\alpha(G-I)\}\|\leq
2|n\overline{q}-\tilde\alpha|+\frac{2\overline{q}(1+p+n\overline
{q}(2-p))}{q+\overline{q}}.\label{serf}
\end{equation}
Even if we choose $\tilde\alpha=n\overline{q}$, the order of
accuracy in
(\ref{serf}) is not better than $n\overline{q}^2$. A similar estimate was
obtained in Theorem 5 of \cite{Wan92}. If we use the terminology of the
book \cite{BHJ92}, we can say that the estimate (\ref{serf}) contains no
`magic' factor.
If we turn to the
papers with `magic' factors, then we have the following results. In
\cite{CM99}, it was proven that if $0\leq p\leq C_0<1$,
then
%e2 ###
%
\begin{equation}\label{CM1}
\|F_n-D_1^n\|\leq
C\max(p_0,\overline{q})\min\biggl(1,\frac{1}{\sqrt{n\overline{q}}} \biggr)+
C\min(\overline{q},n\overline{q}^2)+C\mathrm{e}^{-Cn}.
\end{equation}
The accuracy can be improved, by some asymptotic expansions, to
%e3 ###
%
\begin{eqnarray}\label{CM2}
\bigg\|F_n-D_1^n \biggl(I+p_0\frac{q^2(p-\overline{q})}{(q+\overline{q})^2}I_1\biggr) \bigg\|
&\leq&
C\overline{q}(p+\overline{q})\min\biggl(1,\frac{1}{\sqrt{n\overline{q}}} \biggr)\nonumber
\\[-8pt]\\[-8pt]
&&{}+ C\min(\overline{q},n\overline{q}^2)+C\mathrm{e}^{-Cn}.\nonumber
\end{eqnarray}
Note that in \cite{CM99} formulas (4.5), (4.12) and (4.23) contain
misprints. The parameter $p_0$ is misplaced and should be in the
brackets. If $p\to\tilde p= \mathit{const}$, $\overline{q}\to0$ and
$n\overline{q}\to
\infty$, then the order of accuracy in (\ref{CM1}) is
$\max(\overline{q},(n\overline{q})^{-1/2})$. Also, the order of
accuracy in
(\ref{CM2}) is $\overline{q}$. We can hardly call (\ref{CM2}) the
second order expansion since the improvement of the accuracy was
achieved due to the more precise approximation of the initial
distribution of $\xi_0$ only.

The main idea of signed CP approximations is to leave more than one
factorial cumulant in the exponent. In short, the signed CP measure
has the same structure as the CP measure, but can have negative
Poisson parameters. Such approximations are commonly used in
insurance models and in limit theorems; see \cite
{BC02,Dha94,Hi86,Kor83,Kru86,Pre83,Roo03} and the references therein.
For the
MB distribution in \cite{CM99} the following result is proved. If
%e4 ###
%
\begin{equation}\label{cond0}
p\leq\tilde C<1,\qquad
\frac{\overline{q}}{q+\overline{q}}\leq\frac{1-\tilde C}{30},
\end{equation}
then
%e5 ###
%
\begin{equation}\label{CM4}
\|F_n-D_2^n\|\leq C(p+\overline{q}) \Biggl\{
\min\Biggl(\sqrt{\frac{\overline{q}}{n}},n\overline{q}^2 \Biggr)
+\max(p_0,\overline{q})\min\biggl(1,\frac{1}{\sqrt{n\overline{q}}}\biggr)
+\mathrm{e}^{-Cn} \Biggr\}.
\end{equation}
Note that $\gamma_2<0$ and, therefore, $D_2$ is a signed measure
rather than a distribution. As a rule, signed CP
approximations are more accurate than CP approximations. Indeed, if
$n\to\infty$, then (\ref{CM4}) gives the estimate converging to
zero, even if $p$ and $\overline{q}$ are constants. Also,
(\ref{CM1}) and (\ref{CM2}) are non-trivial even if we only have
$\overline{q}=\mathrm{o}(1)$.

We are unaware of any lower bound estimate for the compound
Poisson approximation to the Markov binomial distribution.
%%%%%%%%%%%%%%%%%%%%%%%%%%%%%%%%%%%%%%%%%%%%%%%%%%%%%%%%%%%%%%%%%%%%%
%%%%%%%%%%%%%%%%%%%%%%%%%%%%%%%%%%%%%%%%%%%%%%%%%%%%%%%%%%%%%%%%%%%%%%
%s3 ###
\section{The main results}
%%%%%%%%%%%%%%%%%%%%%%%%%%%%%%%%%%%%%%%%%%%%%%%%%%%%%%%%%%%%%%%%%%%%%%
%%%%%%%%%%%%%%%%%%%%%%%%%%%%%%%%%%%%%%%%%%%%%%%%%%%%%%%%%%%%%%%%%%%%%%
%s3.1 ###
\subsection{Geometric expansions}
Before formulating our results, it is necessary to
explain the choice of approximating measures. Dobrushin
\cite{Dob53} proved that if $p\to\tilde p$, $n\overline{q}\to\tilde
\lambda$ and $p_0=1$, then the limit distribution for $1+S_{n-1}$
is the convolution $G_1\exp\{\tilde\lambda(G_1-I)\}$,
where $G_1$ is a geometric distribution with parameter $\tilde p$,
that is, $\widehat G_1(t)=(1-\tilde p)\mathrm{e}^{\mathrm
{i}t}/(1-\tilde
p\mathrm{e}^{\mathrm{i}t})$. This suggests that approximation of
$S_n$ for arbitrary $p_0$ should also be based on expansions
in powers of $G-I$.

Let $F$ be concentrated on $\mathbb{Z}$ and have all moments finite.
We can
write formally (i.e., without investigating conditions needed for
the convergence of series)
\[
\widehat F(t)=1+\sum_{j=1}^\infty\frac{\nu_j}{j!}(\mathrm{e}^{\mathrm{i}t}-1)^j=
1+\sum_{m=1}^\infty\frac{\tilde\nu_m}{m!}\bigl(\widehat G(t)-1\bigr)^m.
\]

Here, $\nu_j$ $(j=1,2,\ldots)$ are factorial moments of $F$, and
$\tilde\nu_m$ can be called the \textit{geometric factorial
moments}. Since
\[
\mathrm{e}^{\mathrm{i}t}-1=\frac{q(\widehat G(t)-1)}{1+p(\widehat G(t)-1)},
\]
it is not difficult to establish a relation between $\nu_j$ and
$\tilde\nu_m$:
%e6 ###
%
\begin{equation}\label{g1}
\frac{\tilde\nu_m}{m!}=(-p)^m\sum_{j=1}^m\frac{\nu_j}{j!}
\biggl(-\frac {q}{p} \biggr)^j\pmatrix{m-1 \cr m-j},\qquad m=1,2,\ldots.
\end{equation}
Similar relations hold for factorial cumulants and geometric
factorial cumulants. For the MB distribution, we have
%e7 ###
%
\begin{equation}
\tilde\nu_1=q\nu_1= q\mathrm{E}
S_n=n\gamma_1+\varkappa_1+\varkappa_2-(\varkappa_1+\varkappa
_2)(p-\overline{q})^n.
\label{g2}
\end{equation}
(For the formula of the mean, see \cite{MBB}.) Since we will
assume $p$ and $\overline{q}$ to be small, the last summand in
(\ref{g2}) will be neglected. As it turns out, $F_n$ is close to
some convolution $W_1\Lambda_1^n$; see (\ref{a10}) below. We use
(\ref{g1}) for choosing the approximating measure for $W_1$. The cumulant
analog of (\ref{g1}) is used for $\Lambda_1^n$.

%s3.2 ###
\subsection{Compound Poisson approximation}

In this paper, we usually assume that
%e8 ###
%
\begin{equation}\label{cond1}
p\leq\tfrac{1}{2},\qquad \overline{q}\leq\tfrac{1}{30}.
\end{equation}
The size of the absolute constants is determined by the
method of proof. We expand $F_n$ as a series of convolutions of
measures. The remainder term is usually estimated by a series
containing powers of $\overline{q}+p$. If the sum $\overline{q}+p$ is
sufficiently small, the series converges. Thus, although we have some
freedom in the choice of magnitude of $p$ and $\overline{q}$, the sum
$p+\overline{q}$ must be small. The choice of condition (\ref{cond1}) is
determined by the fact that the CP limit occurs when
$n\overline{q}\to\tilde\lambda$; see, for example, Table 1 in \cite{Dob53}.
Therefore, if we expect the CP approximation to be accurate, then
$\overline{q}$ should be small. On the other hand, we have included the
case $\overline{q}=\mathit{constant}$, that is, the case which is usually associated
with the normal approximation.
We choose the
assumption $p\leq1/2$ instead of (\ref{cond0}), in order to make
our proofs clearer.
%Note, by the way, that for $\tilde C=1/2$,
%(\ref{cond0}) gives $\qubar\leq1/59$ in contrast to our
%$\qubar\leq1/30$.

%THEOREMA
%%%%%%%%%%%%%%%%%%%%%%%%%%%%%%%%%%%%%%%%%%%%%%%%%%%%%%%%%%%%%%%%%%%%%%
%
%t3.1
\begin{theorem} \label{T1} Let $p\leq1/2$. We then have
%e9 ###
%
\begin{equation}\label{CP1}
\|F_n-HD_1^\lambda\|\leq
C\overline{q}(p+\overline{q})\min\biggl(1,\frac{1}{\sqrt{n\overline
{q}}} \biggr)+
C\min(\overline{q},n\overline{q}^2)+C(p+\overline{q})\mathrm{e}^{-C_1n}.
\end{equation}
If, in addition, $\overline{q}\leq1/30$, then
%e11 ###
%e10 ###
%
\begin{eqnarray}\label{CP2}
\|F_n-HD_1^\lambda\|_\infty
&\leq&
C\overline{q}(p+\overline{q})\min\biggl(1,\frac{1}{{n\overline{q}}} \biggr)+
C\min\Biggl(\sqrt{\frac{\overline{q}}{n}},n\overline{q}^2 \Biggr)\nonumber\\[-8pt]\\[-8pt]
&&{}+C(p+\overline{q})\mathrm{e}^{-C_1n},\nonumber\label{CP3}
\\
\|F_n-HD_1^\lambda\|_{\mathrm{W}}&\leq& C\overline{q}(p+\overline{q})+
C\min\bigl(\overline{q}\sqrt{n\overline{q}},n\overline
{q}^2\bigr)+C(p+\overline{q})\mathrm{e}^{-C_1n}.
\end{eqnarray}
\end{theorem}

%c3.1
\begin{corollary} \label{cor1} If (\ref{cond1}) is satisfied and
$n\overline{q}\geq1$, then
\[%\label{CP4}
\|F_n-HD_1^n\|\leq C\overline{q},\qquad
\|F_n-HD_1^n\|_\infty\leq C\sqrt{\frac{\overline{q}}{n}},\qquad
\|F_n-HD_1^n\|_{\mathrm{W}}\leq C\overline{q}\sqrt{n\overline{q}}.
\]
\end{corollary}

%c3.2
\begin{corollary} \label{cor2} If (\ref{cond1}) is satisfied, then
\[
\|F_n-HH_1^\lambda\|\leq C(\overline{q}+p\mathrm{e}^{-C_1n}).\label{CP6}
\]
\end{corollary}

%r3.1
\begin{Remarks}\vspace*{-6pt}
\begin{enumerate}[(iii)]
\item[(i)] $H$ is compound
Poisson distribution; see Lemma~\ref{trys} below. If
$\overline{q}\leq p$, then $H_1^\lambda$ is also a CP distribution.
Thus, we see that there exist quite different forms of CP
approximations with similar orders of accuracy.
\item[(ii)] The estimate (\ref{CP1}) is slightly better than
(\ref{CM2}) for $p,\overline{q}\leq\exp\{-C_1n\}$, and more
accurate than (\ref{CM1}) for $\overline{q}\leq1/\sqrt{n\overline{q}}$
and $p_0\geq C$.
\item[(iii)] For the closeness of $F_n$ and $H_1D_1^\lambda$, it
suffices to assume $\overline{q}\to0$, in considerable contrast to
$n\overline{q}\to
\tilde\lambda$, the latter being needed for the convergence to the
limit CP
law.
\item[(iv)] We can write
$\min(\overline{q},n\overline{q}^2)=n\overline{q}^2\min
(1,(n\overline{q})^{-1})$. The last
factor, in terms of \cite{BHJ92}, page~5, can be called the `magic' factor.
\end{enumerate}
\end{Remarks}

The accuracy of approximation can be improved by the second order
expansion.
%t3.2
\begin{theorem}\label{T2} If $p\leq1/2$, then
\begin{eqnarray*}
\big\|F_n-HD_1^\lambda\bigl(I+n\gamma_2(G-I)^2\bigr)\big\|
\leq C \biggl\{\overline{q}^2+p\overline{q}\min\biggl(1,\frac{1}{\sqrt{n\overline{q}}} \biggr)+(p+\overline{q})\mathrm{e}^{-C_1n} \biggr\}.
\label{CP7}
\end{eqnarray*}
If, in addition, $\overline{q}\leq1/30$, then
\begin{eqnarray*}
&&\big\|F_n-HD_1^\lambda\bigl(I+n\gamma_2(G-I)^2\bigr)\big\|_\infty\nonumber
\\
&&\quad\leq
C \biggl\{\overline{q}^2\min\biggl(1,\frac{1}{\sqrt{n\overline{q}}}\biggr)
+p\overline{q}\min\biggl(1,\frac{1}{n\overline{q}} \biggr)+(p+\overline
{q})\mathrm{e}^{-C_1n} \biggr\},
\label{CP8}\\
&&\big\|F_n-HD_1^\lambda\bigl(I+n\gamma_2(G-I)^2\bigr)\big\|_{\mathrm{W}}\nonumber\\
&&\quad\leq
C\bigl\{\overline{q}^2\max\bigl(1,\sqrt{n\overline{q}}\bigr)+p\overline{q}+(p+\overline{q})\mathrm{e}^{-C_1n}\bigr\}.
\label{CP9}
\end{eqnarray*}
\end{theorem}

Note that the last estimate contains $\max(1, \sqrt{n\overline{q}})$,
reflecting the fact that the estimates for the Wasserstein
distance are less accurate than the ones for the total variation
norm. It is even more evident when $n\overline{q}\geq1$.
%
%c3.3
\begin{corollary} If (\ref{cond1}) is satisfied and
$n\overline{q}\geq1$, then the estimates in Theorem~\ref{T2} are
\[
C\overline{q}\biggl(\overline{q}+\frac{p}{\sqrt{n\overline{q}}} \biggr),\qquad
C\sqrt{\frac{\overline{q}}{n}} \biggl(\overline{q}+\frac{p}{\sqrt{n\overline{q}}} \biggr),\qquad
C\overline{q}\sqrt{n\overline{q}} \biggl(\overline{q}+\frac{p}{\sqrt{n\overline{q}}} \biggr),
\]
respectively.
\end{corollary}

We see that, in general, even the second order estimates in total
variation are only meaningful for $\overline{q}=\mathrm{o}(1)$.

%s3.3 ###
\subsection{Signed compound Poisson approximations}

The choice of a signed CP approximation, in general, means
that the first term of the asymptotic expansion, unlike Theorem
\ref{T2}, is in the exponent.
%t3.3
\begin{theorem}\label{T3} If condition (\ref{cond1}) is satisfied,
then
%e12 ###
%
\begin{eqnarray}
\|F_n-H\exp\{\varkappa_1(G-I)\}D_2^n\|&\leq&
C(p+\overline{q}) \Biggl\{\min\Biggl(\overline{q},\sqrt{\frac{\overline
{q}}{n}} \Biggr)+\mathrm{e}^{-C_1n} \Biggr\},
\label{SCP1}\\
\|F_n-H\exp\{\varkappa_1(G-I)\}D_2^n\|_\infty&\leq&
C(p+\overline{q}) \biggl\{\min\biggl(\overline{q},\frac{1}{n} \biggr)+\mathrm
{e}^{-C_1n} \biggr\},
\nonumber
\\
\|F_n-H\exp\{\varkappa_1(G-I)\}D_2^n\|_{\mathrm{W}}&\leq&
C(p+\overline{q})\{\overline{q}+\mathrm{e}^{-C_1n}\}. \nonumber
\end{eqnarray}
\end{theorem}

Note that for $n\overline{q}\leq1$ and $p_0=\mathit{constant}$, (\ref{SCP1})
is more accurate that (\ref{CM4}). More importantly, when
$p=\mathit{constant}$ and $\overline{q}=\mathit{constant}$, the estimate (\ref{SCP1})
is of
order $\mathrm{O}(n^{-1/2})$. In this sense, the signed CP approximation is
comparable to the normal one and, moreover, it holds in the total
variation metric. Meanwhile, for discrete distributions, the
normal approximation holds in the uniform metric only. Just as
in the CP case, the second order expansions can be used.
%
%t3.4
\begin{theorem}\label{T4} If (\ref{cond1}) holds, then
%e15 ###
%e14 ###
%e13 ###
%
\begin{eqnarray}\label{SCP4}
&&\big\|F_n-H\exp\{\varkappa_1(G-I)\}D_2^n\bigl(I+n\gamma_3(G-I)^3\bigr)\big\|
\nonumber\\[-8pt]\\[-8pt]
&&\quad\leq
C(p+\overline{q}) \biggl\{\min\biggl(\overline{q},\frac{1}{n} \biggr)+\mathrm{e}^{-C_1n} \biggr\},\nonumber\\
&&\big\|F_n-H\exp\{\varkappa_1(G-I)\}D_2^n\bigl(I+n\gamma_3(G-I)^3\bigr)\big\|_\infty\label{SCP5}
\nonumber\\[-8pt]\\[-8pt]
&&\quad\leq
C(p+\overline{q}) \biggl\{\min\biggl(\overline{q},\frac{1}{n\sqrt{n\overline
{q}}} \biggr)+\mathrm{e}^{-C_1n} \biggr\}
%%%%%%%%%%%%%%%%%%%%%%%%%%%%%%%%%%%%%%%%%%%%%%%%%%%%%%%%%%%%% \m
,\nonumber
\\
&&\big\|F_n-H\exp\{\varkappa_1(G-I)\}D_2^n\bigl(I+n\gamma_3(G-I)^3\bigr)\big\|_{\mathrm{W}}\label{SCP6}
\nonumber\\[-8pt]\\[-8pt]
&&\quad\leq
C(p+\overline{q}) \Biggl\{\min\Biggl(\overline{q},\sqrt{\frac{\overline
{q}}{n}} \Biggr)+\mathrm{e}^{-C_1n} \Biggr\}.\nonumber
\end{eqnarray}
\end{theorem}

%c3.4
\begin{corollary} Let $n\overline{q}\geq1$. The estimates (\ref{SCP4})--(\ref{SCP6}) are then at
least of order
\[
\frac{C(p+\overline{q})}{n},\qquad
\frac{C(p+\overline{q})}{n\sqrt{n\overline{q}}},\qquad
\frac{C(p+\overline{q})\sqrt{\overline{q}}}{\sqrt{n}},
\]
respectively.
\end{corollary}

In Theorem~\ref{T4}, only a part of the asymptotic expansion is in the
exponent. Therefore, the following question naturally arises. Is it
possible to find a signed CP measure which, up to a constant,
provides the same accuracy as in Theorem~\ref{T4}? As it follows
from the following result, such a measure indeed exists.
%t3.5
\begin{theorem}\label{T5} If (\ref{cond1}) holds, then
\begin{eqnarray*}
\|F_n-H\exp\{\varkappa_1(G-I)\}D_3^n\|
&\leq&
C(p+\overline{q}) \biggl\{\min\biggl(\overline{q},\frac{1}{n} \biggr)+\mathrm{e}^{-C_1n} \biggr\}, \label{SCP7}\\
\|F_n-H\exp\{\varkappa_1(G-I)\}D_3^n\|_\infty
&\leq&
C(p+\overline{q}) \biggl\{\min\biggl(\overline{q},\frac{1}{n\sqrt{n\overline
{q}}} \biggr)+\mathrm{e}^{-C_1n} \biggr\},\label{SCP8}\\
\|F_n-H\exp\{\varkappa_1(G-I)\}D_3^n\|_{\mathrm{W}}
&\leq&
C(p+\overline{q}) \Biggl\{\min\Biggl(\overline{q},\sqrt{\frac{\overline
{q}}{n}} \Biggr)+\mathrm{e}^{-C_1n} \Biggr\}.\label{SCP9}
\end{eqnarray*}
\end{theorem}

%%%%%%%%%%%%%%%%%%%%%%%%%%%%%%%%%%%%%%%%%%%%%%%%%%%%%%%%%%%%%%%%%%%%%%%%%%%%%%%%%%%%%%%%%%%%%%%%%%%%%
%%%%%%%%%%%%%%%%%%%%%%%%%%%%%%%%%%%%%%%%%%%%%%%%%%%%%%%%%%%%%%%%%%%%%%%%%%%%%%%%%%%%%%%%%%%%%%%%%%%%%
%%%%%%%%%%%%%%%%%%%%%%%%%%%%%%%%%%%%%%%%%%%%%%%%%%%%%%%%%%%%%%%%%%%%%%%%%%%%%%%%%%%%%%%%%%%%%%%%%%%%%
%s3.4 ###
\subsection{Lower bound estimates}
In this section, we show that in some cases, the estimates in
Theorems~\ref{T1} and~\ref{T3} are of the correct order. We
concentrate our attention on the case $n\overline{q}\geq1$.
%t3.6
\begin{theorem}\label{T6} Let condition (\ref{cond1}) be satisfied
and let $n\overline{q}\geq1$. Then, for some absolute constants
$C_2$ and $C_3$,
%e18 ###
%e17 ###
%e16 ###
%
\begin{eqnarray}
\|F_n-HD_1^\lambda\|&\geq&
C_2\overline{q}\biggl(1-C_3 \biggl(\overline{q}+\frac{p}{\sqrt{n\overline{q}}}\biggr) \biggr),\label{LB1}\\
\|F_n-HD_1^\lambda\|_\infty&\geq&
C_2\sqrt{\frac{\overline{q}}{n}} \biggl(1-C_3 \biggl(\overline{q}+\frac
{p}{\sqrt{n\overline{q}}} \biggr) \biggr),\label{LB2}\\
\|F_n-HD_1^\lambda\|_{\mathrm{W}}&\geq&
C_2\overline{q}\sqrt{n\overline{q}} \biggl(1-C_3 \biggl(\overline{q}+\frac
{p}{\sqrt{n\overline{q}}} \biggr) \biggr).\label{LB3}
\end{eqnarray}
\end{theorem}

It is obvious that estimates (\ref{LB1})--(\ref{LB3}) are
non-trivial only when the expression in the brackets is positive. Let
$p\leq1/2$, $n\overline{q}\to\infty$ and $\overline{q}\to0$. Combining
Theorems~\ref{T1} and~\ref{T6}, for sufficiently large $n$, we
obtain
\begin{eqnarray*}
C_4\overline{q}&\leq&\|F_n-HD_1^\lambda\|\leq C_5\overline{q}, \\
C_4\sqrt{\frac{\overline{q}}{n}}&\leq&\|F_n-HD_1^\lambda\|_\infty\leq C_5\sqrt{\frac{\overline{q}}{n}},\\
C_4\overline{q}\sqrt{n\overline{q}}&\leq&\|F_n-HD_1^\lambda\|_{\mathrm{W}}\leq C_5\overline{q}\sqrt{n\overline{q}}.
\end{eqnarray*}
Of course, the last estimate, as well as the one in (\ref{LB3}),
is of interest only if $\overline{q}\sqrt{n\overline{q}}\to0$. Similar
results can be obtained for the signed CP approximations.
%
%t3.7
\begin{theorem}\label{T7} Let condition (\ref{cond1}) be satisfied
and let $n\overline{q}\geq1$. Then, for some absolute constants
$C_6$ and $C_7$,
%e21 ###
%e20 ###
%e19 ###
%
\begin{eqnarray}
\|F_n-H\exp\{\varkappa_1(G-I)\}D_2^n\|
&\geq& C_6\sqrt{\frac{\overline{q}}{n}} \biggl(|\tilde\gamma_3|-C_8\frac{p+\overline{q}}{\sqrt{n\overline{q}}} \biggr),\label{LB4}\\
\|F_n-H\exp\{\varkappa_1(G-I)\}D_2^n\|_\infty
&\geq&\frac{C_6}{n} \biggl(|\tilde\gamma_3|-C_8\frac{p+\overline{q}}{\sqrt{n\overline{q}}} \biggr),\label{LB5}\\
\|F_n-H\exp\{\varkappa_1(G-I)\}D_2^n\|_{\mathrm{W}}
&\geq& C_6\overline{q}\biggl(|\tilde\gamma_3|-C_8\frac{p+\overline{q}}{\sqrt{n\overline{q}}} \biggr).\label{LB6}
\end{eqnarray}
\end{theorem}

Let $n\overline{q}\to\infty$ as $\overline{q}\to0$ and $p\to
\tilde p$. Also,
assume that $n$ is sufficiently large so that the right-hand estimates
of (\ref{LB4})--(\ref{LB6}) are positive. We then have
\begin{eqnarray*}
C_8\sqrt{\frac{\overline{q}}{n}}&\leq&\|F_n-H\exp\{\varkappa_1(G-I)\}D_2^n\|\leq C_9\sqrt{\frac{\overline{q}}{n}},\\
\frac{C_8}{n}&\leq&\|F_n-H\exp\{\varkappa_1(G-I)\}D_2^n\|_\infty\leq\frac{C_9}{n},\\
C_8\overline{q}&\leq&\|F_n-H\exp\{\varkappa_1(G-I)\}D_2^n\|_{\mathrm{W}}\leq C_9\overline{q}.
\end{eqnarray*}

%s3.5 ###
\subsection{Asymptotically sharp constants}

In the previous section, we proved that upper and lower bound
estimates are of the same order, provided that $n\overline{q}$ is
large and
$\overline{q}$ is small. As it turns out, if, in addition, $p$ is small,
then it is possible to obtain asymptotically sharp constants.
%t3.8
\begin{theorem}\label{T8} Let $p\leq1/4$, $\overline{q}\leq
1/30$ and $n\overline{q}\geq1$. Then
%e24 ###
%e23 ###
%e22 ###
%
\begin{eqnarray}
\big\vert\|F_n-HD_1^\lambda\|-A_{11} \big\vert
&\leq&C\overline{q}\biggl(p+\overline{q}+\frac{1}{\sqrt{n\overline{q}}}\biggr),\label{ASH1}\\
\big\vert\|F_n-HD_1^\lambda\|_\infty-A_{12} \big\vert
&\leq&C\sqrt{\frac{\overline{q}}{n}} \biggl(p+\overline{q}+\frac{1}{\sqrt{n\overline{q}}}\biggr)
,\label{ASH2}\\
\big\vert\|F_n-HD_1^\lambda\|_{\mathrm{W}}-A_{13} \big\vert
&\leq&C\overline{q}\sqrt{n\overline{q}} \biggl(p+\overline{q}+\frac{1}{\sqrt{n\overline{q}}}\biggr)
,\label{ASH3}
\end{eqnarray}
where
\[
A_{11}=\frac{4|\gamma_2|}{\gamma_1 q\sqrt{2{\curpi}\mathrm{e}}},\qquad
A_{12}=\frac{|\gamma_2|}{\gamma_1\sqrt{\gamma_1}\sqrt{2{\curpi}n q}},\qquad
A_{13}=\frac{|\gamma_2|\sqrt{2n}}{q\sqrt{\gamma_1{\curpi}q}}.
\]
\end{theorem}

As a consequence of~(\ref{ASH1}), we note that if
$p\to0$, $\overline{q}\to0$ and
$n\overline{q}\to\infty$, then
\[
\|F_n-HD_1^\lambda\|\sim\frac{6\overline{q}}{\sqrt{2{\curpi}
\mathrm{e}}}.
\]
Similar relations can be obtained for the local and Wasserstein
norms as well.

%%%%%%%%%%%%%%%%%%%%%%%%%%%%%%%%%%%%%%%%%%%%%%%%%%%%%%%%%%%%%%%%%%%%%%

%s4 ###
\section{Applications of Markov binomial models}

In this section, we discuss some areas where the results of our
paper can be applied.
\begin{longlist}[(iii)]
\item[(i)] \textit{Aggregate claim distribution in the
individual model.} Consider a portfolio of $n$ risks. Each risk
produces a positive claim amount during a certain reference
period. The aggregate claim of the portfolio is then
\[
S^{\mathrm{ind}}=X_1+X_2+\cdots+X_n.
\]
It is usually assumed that all $X_j$ are independent. However, the
independence of claims does not always reflect reality. For
example, an accident involving a tourist group, life insurance for a husband
and wife or pensions for workers of the same company are likely
to produce dependent risks. For discussion of the dependence of
risks and further examples, see \cite{GDH96,DHG97} and
\cite{Ribas03}.

Compound Poisson and signed compound Poisson approximations in the
independent case of an individual model have been quite thoroughly
investigated; see, for example, \cite{Ge84,Hi86}. On the
other hand, there are only a few results for the total variation
metric for dependent risks. Dhaene and Goovaerts \cite{DHG97}
investigated a similar model (although not explicitly Markovian) under
an assumption which, in our notation, is equivalent to $\overline{q}_m=0$.
However, under such an assumption, one cannot expect the limiting law
to be compound Poisson. Therefore, we have excluded this peripheral case
from this paper, assuming $\overline{q}$ to be small, but not identically
zero.
In \cite{GDH96}, Poisson approximation in the general setting of
dependent risks was discussed.
However, in our case, their result is not applicable since for small
$\overline{q}$, the distribution of the approximated sum is not close
to the Poisson
distribution, but rather to
the compound Poisson law.

Let us assume that aggregated claim amount $S^{\mathrm{ind}}$ of the
portfolio consists of $N$ independent groups of risks. We assume
a homogeneous model for each group of risks with Markovian
dependence. Let each risk have a two-point distribution. More
precisely, let
\[
S^{\mathrm{ind}}=\sum_{m=1}^N\sum_{j=1}^{n_m}X_{j}^m.
\]
Here, $X_j^m$ and $X_k^l$ are independent if $m\neq l$. We assume
that each risk of the $m$th group can produce a claim of size
$a_m$. Moreover, the dependence of risks of the same group is
Markovian: $\mathrm{P}(X_1^m=a_m)=\overline{q}_m$, $\mathrm
{P}(X_1^m=0)=\overline{p}_m$
and
\begin{eqnarray*}
%& \Prob(X_1^m=a_m)=\qubar_m,
 \mathrm{P}(X_j^m=a_m | X_{j-1}^m=0)&=&\overline{q}_m,\qquad
\mathrm{P}(X_j^m=0 | X_{j-1}^m=0)=\overline{p}_m,\\
 \mathrm{P}(X_j^m=a_m | X_{j-1}^m=a_m)&=&p_m<1/2,\qquad
\mathrm{P}(X_j^m=0 | X_{j-1}^m=a_m)=q_m,\\[-28pt]
\end{eqnarray*}
\begin{eqnarray*}
p_m+q_m=\overline{q}_m+\overline{p}_m=1,\qquad
p_m,\overline{q}_m\in(0,1),\qquad m=1,2,\ldots,N, j=2,\ldots,n_m.
\end{eqnarray*}
The results of the previous sections
can now easily be applied.
We illustrate this with just one example.
Let us define a compound Poisson variable in the following way:
\[
S^{\mathrm{cp}}=\sum_{m=1}^N a_m\sum_{j=0}^{N_m}Y_{jm}.
\]
Here, $Y_{jm}$ are i.i.d.~geometric random variables,
$\mathrm{P}(Y_{jm}=k)=q_mp_m^{k-1}$, $k=1,2,\ldots,$ and $N_m$ is a
Poisson random variable with parameter
$n_mq_m\overline{q}_m/(q_m+\overline{q}_m)$. The random variables $N_m$,
$m=1,2,\ldots,N$, are independent and also do not depend on
$Y_{jm}$.
Denote the distributions of $S^{\mathrm{ind}}$ and $S^{\mathrm{cp}}$ by $F^{\mathrm{ind}}$ and
$F^{\mathrm{cp}}$, respectively. The characteristic function of $F^{\mathrm{cp}}$
is then given by
\[
\widehat
F^{\mathrm{cp}}(t)=\exp\Biggl\{\sum_{m=1}^N\frac{n_mq_m\overline{q}_m(\mathrm
{e}^{\mathrm{i}t a_m }-1)}{(q_m+\overline{q}_m)(1-p_m\mathrm
{e}^{\mathrm{i}t a_m })} \Biggr\}.
\]
Also, we have the following estimate of approximation:
%e25 ###
%
\begin{eqnarray}\label{ACP1}
\|F^{\mathrm{ind}}-F^{\mathrm{cp}}\|&\leq&
C\sum_{m=1}^N [\overline{q}_m(p_m+\overline{q}_m)\min
(1,(n_m\overline{q}_m)^{-1/2})\nonumber\\[-8pt]\\[-8pt]
&&{}+
\min(\overline{q}_m,n_m\overline{q}_m^2)+(p_m+\overline
{q}_m)\mathrm{e}^{-C_1n_m} ].\nonumber
\end{eqnarray}
Note that the approximation is closer if all $\overline{q}_m$ are
small.
%In comparison to model from [Dh G1997] with $\qubar=0$,
%this is a mild restriction.

For the proof of (\ref{ACP1}), one should use the triangle
inequality, thus reducing the problem to $N$ estimates of Markov
binomial distributions concentrated on $0,a_m,2a_m,\ldots.$ The
total variation metric is invariant with respect to norming.
Therefore, without loss of generality, one can switch to integer
numbers and apply (\ref{CP1}) $N$ times with $p_0=0$.

It is obvious that the second order estimates and estimates in
Wasserstein metric can be obtained in a similar way.

\item[(ii)] \textit{System failure models.} The Markov binomial
distribution
naturally arises in weather and stock market trends. It is also a
natural model for system failure situations. As an example, we
present one model from Sahinoglu \cite{Sog90}, who considered
an electric power supply system with operating and non-operating
states throughout a year-long period of operation, discretized in
hours. Let $M_i$ be the margin values at hourly steps, that is,
\[
M_i=\mathit{TPG}-X-L_i,
\]
where $\mathit{TPG}$ denotes total power generation, $L_i$ denotes power
demand (hourly peak load forecast) and $X$ denotes unplanned
forced outages. Let $Y_i$ be an indicator of $\{M_i<0\}.$ Then
$S=Y_1+Y_2+\cdots+Y_n$ represents cumulated hours of negative-margin
hours, that is, the unavailability of power at the $n$th hour.
It is natural to assume that $S$ has a Markov binomial
distribution. Notably, the Markov chain $Y_1,Y_2,\ldots$ is
non-stationary. Therefore, many known results about the compound
Poisson approximation cannot be applied directly. Also, the
results of our paper relax the assumptions on transition
probabilities from Sahinoglu's model and give estimates of the
accuracy of approximations. Further, as shown in Sahinoglu \cite{Sog90},
page 49, the probabilities of the compounding geometric law under
certain assumptions can be viewed as probabilities for the number of
trials required to repair the system.

\item[(iii)] \textit{Industrial applications: sampling plans.}
A basic assumption in standard acceptance plans for attributes is
that the characteristics of items in the lots are i.i.d.~Bernoulli
variables. Recently, however, the focus has been on monitoring the
ongoing production process by inspecting the items sequentially.
In such cases, the quality levels of successive items are statistically
dependent and it has been found in practice that the Markov-dependent
model is a very useful one; see \cite{Nel93}.
Indeed, Bhat \textit{et al.}~\cite{BhL90}
modified the standard acceptance sampling plans and proposed
sequential single sampling plans for monitoring Markov-dependent
production processes. Vellaisamy and Sankar \cite{Vel01} proposed
optimal systematic sampling plans for Markov-dependent processes.
We will outline some possible new research directions in this
field.
\end{longlist}
%s5 ###
\section{Auxiliary results}

We now introduce further notation:
%e28 ###
%e27 ###
%e26 ###
%
\begin{eqnarray}
a_1&=&\gamma_1,\qquad
a_2=\gamma_2+\frac{a_1^2}{2},\qquad
a_3=\gamma_3+a_1a_2-\frac{a_1^3}{3},\label{a1}\\
Y&=&G-I,\qquad
B=\sum_{j=0}^\infty(pI_1-\overline{q}I)^j,\qquad
K=\sum_{j=0}^\infty(pI_1-\overline{q}I-2\gamma_1Y)^j,\label{a23}\\
L&=&\frac{4\overline{q}^2}{(q+\overline{q})^2}Y^2
[q^2I+p(q+\overline{q})(I-pI_1)
]K^2.\label{a4}
\end{eqnarray}
In the following two lemmas, $C(k)$ denotes an absolute positive
constant depending on $k$. Throughout this paper, we set $0^0=1$.

%%%%%%%%%%%%%%%%%%%%%%%%%%%%%%%%%%%%%%%%%%%%%%%%%%%%%%%%%%%%%%%%%%%%%%
%
%l5.1
\begin{lemma} \label{smooth}
Let $t>0$, $k \in\{0, 1, \ldots\}$ and $0<p<1$. Also, let $M$
be a finite (signed) measure concentrated at $\mathbb{Z}$. Then, for
$Y$ defined in (\ref{a23}),
%e29 ###
%
\begin{equation}\label{a5}
\|Y^2\mathrm{e}^{tY}\|\leq\frac{3}{t\mathrm{e}},\qquad
\|Y^k\mathrm{e}^{tY}\|\leq\biggl(\frac{2k}{t\mathrm{e}} \biggr)^{k/2}.
\end{equation}
If $p\leq1/2$, then
%e30 ###
%
\begin{equation}\label{a67}
\|Y^k\mathrm{e}^{tY}\|_\infty\leq\frac{C(k)}{t^{(k+1)/2}},\qquad
\|YM\|_{\mathrm{W}}\geq\frac{2}{3}\|M\|,\qquad
\|Y M\|\geq\frac{2}{3}\|(I_1-I)M\|.
\end{equation}
\end{lemma}

%%%%%%%%%%%%%%%%%%%%%%%%%%%%%%%%%%%%%%%%%%%%%%%%%%%%%%%%%%%%%%%%%%%%%%
%
\begin{pf}
The estimates in (\ref{a5}) follow from the properties of the
total variation norm and results in \cite{Roo01} and \cite{DP88}.
The first estimate in (\ref{a67}) is a consequence of the inversion
formula and the following inequalities:
\[
\operatorname{Re}\widehat Y(t)\leq-\frac{2}{1+p}\sin^2\frac{t}{2},\qquad
|\widehat Y(t)|\leq\frac{2}{q}\bigg\vert\sin\frac{t}{2} \bigg\vert.
\]
Here, $\operatorname{Re}\{\cdot\}$ means the real part of the complex number. In
view of the relation between total variation and Wasserstein norms
(see the \hyperref[sec1]{Introduction}), we get
\begin{eqnarray*}
\|MY\|_{\mathrm{W}}&=&\Bigg\|(I_1-I)M\sum_{j=0}^\infty p^jI_j\Bigg\|_{\mathrm{W}}=\Bigg\|M\sum_{j=0}^\infty p^jI_j\Bigg\|,\\
\|M(I_1-I)\|&=&\|MY(I-pI_1)\|\leq\|MY\|(1+p),\\
\|M\|&=&\Bigg\|M\sum_{j=0}^\infty p^jI_j(I-pI_1)\Bigg\|\leq\Bigg\|M\sum_{j=0}^\infty p^jI_j\Bigg\|(1+p).
\end{eqnarray*}
The results in (\ref{a67}) now follow easily.
\end{pf}

For our asymptotically sharp results, we need the following lemma.
Set
\begin{eqnarray*}
&&\varphi_k(x)=\frac{1}{\sqrt{2{\curpi}}}\frac{{\mathrm{d}}^k}{{\mathrm{d}}x^k} \mathrm{e}^{-x^2/2},\qquad
\|\varphi_k\|_1=\int_\mathbb{R}|\varphi_k(x)| \,{\mathrm{d}}x,\qquad
\|\varphi_k\|_\infty=\sup_{x\in\mathbb{R}}|\varphi_k(x)|\nonumber
\\
&&\quad(k=0,1,\ldots).
\end{eqnarray*}
%
%%%%%%%%%%%%%%%%%%%%%%%%%%%%%%%%%%%%%%%%%%%%%%%%%%%%%%%%%%%%%%%%%%%%%%
%
%l5.2
\begin{lemma}\label{sharpC}
Let $t>0$ and $k=0,1,2,\ldots.$ We then have
\begin{eqnarray*}
\bigg\vert\big\|(I_1-I)^k\mathrm{e}^{t(I_1-I)}\big\|-\frac{\|\varphi_{k}\|_1}{t^{k/2}} \bigg\vert &\leq&\frac{C(k)}{t^{(k+1)/2}},\\
\bigg\vert\big\|(I_1-I)^k\mathrm{e}^{t(I_1-I)}\big\|_\infty-\frac{\|\varphi_{k}\|_\infty}{ t^{(k+1)/2}} \bigg\vert&\leq&\frac{C(k)}{t^{k/2+1}},\\
\bigg\vert\big\|(I_1-I)^k\mathrm{e}^{t(I_1-I)}\big\|_{\mathrm{W}}-\frac{\|\varphi_{k-1}\|_1}{t^{(k-1)/2}} \bigg\vert&\leq&\frac{C(k)}{t^{k/2}}\qquad
(k\neq0).
\end{eqnarray*}
\end{lemma}

The proof follows from a more general Proposition 4 in
\cite{Roo03}.

%l5.3
\begin{lemma}\label{trys} If $N>0$ and $0<\alpha\leq p<1$, then
$(I+\alpha Y)^N$ is a CP distribution.
\end{lemma}

\begin{pf}
Note that
\[
(I+\alpha Y)^N=\exp\{-N\ln(1-\alpha)(F-I) \}.
\]
Here, $F$ is a distribution concentrated on $\{1,2,3,\ldots\}$ with
\[
F\{j\}=-\frac{1}{\ln(1-\alpha)}\frac{1}{j} \biggl(p^j- \biggl(\frac{p-\alpha}{1-\alpha} \biggr)^j \biggr).
\]
The last relation obviously completes the proof.
\end{pf}

Before we proceed to our main lemma, we need some additional facts
about $F_n$. Similarly to \cite{CM99} (see also \cite{MBB}), it
is possible to check that under assumption (\ref{cond1}), we have
%e31 ###
%
\begin{equation}\label{a9}
\widehat{F}_n(t)=\widehat{\Lambda}_1^n(t) \widehat{W}_1(t)
+\widehat{\Lambda}_2^n(t) \widehat{W}_2(t),
\end{equation}
where
\begin{eqnarray*}
\widehat{\Lambda}_{1,2}(t)&=&\frac{p\mathrm{e}^{\mathrm{i}t}+\overline{p}\pm
\widehat{D}^{1/2}(t)}{2},\label{wL12}\\
\widehat{W}_{1,2}(t)&=&\frac{p_0}{2} \biggl(1\pm\frac{q+\overline{q}+p(\mathrm{e}^{\mathrm{i}t}-1)}{\widehat{D}^{1/2}(t)} \biggr)
+\frac{1-p_0}{2} \biggl(1\pm\frac{q+\overline{q}+(2\overline{q}-p)(\mathrm{e}^{\mathrm{i}t}-1)}{
\widehat{D}^{1/2}(t)} \biggr),\label{wW12}\\
\widehat{D}(t)&=&(p\mathrm{e}^{\mathrm{i}t}+\overline{p})^2+4\mathrm{e}^{\mathrm{i}t}(\overline{q}-p).
\end{eqnarray*}
This allows us to write $F_n$ as
%e32 ###
%
\begin{equation}
F_n=\Lambda_1^nW_1+\Lambda_2^nW_2 \label{a10}
\end{equation}
and to express $\Lambda_{1,2}$ and $W_{1,2}$ as the following
series:
%%%%%%%%%%%%%%%%%%%%%%%%%%%%%%%%%%%%%%%%%%%%%%%%%%%%%%%%%%%%%%%%%%%%%%
%%%%%%%%%%%%%%%%%%%%%%%%%%%%%%%%%%%%%%%%%%%%%%%%%%%%%%%%%%%%%%%%%%%%%%
%e35 ###
%e34 ###
%e33 ###
%
\begin{eqnarray}
\Lambda_1&=&I+a_1Y+\frac{1}{2} \{(1+\overline{q})I-pI_1+2a_1Y\}\sum
_{j=1}^\infty\pmatrix{1/2 \cr j}
(-1)^jL^j, \label{a11}\\
\Lambda_2&=&pI_1-\overline{q}I+(I-\Lambda_1),\label{a12}\\
\label{a13}
W_{1,2}&=&\frac{1}{2} \Biggl\{I\pm[(q+\overline{q})I+p(I_1-I)]K\sum_{j=0}^\infty\pmatrix{-1/2 \cr j}
(-1)^jL^j \Biggr\}\nonumber\\[-8pt]\\[-8pt]
&&{}\pm(1-p_0)(\overline{q}-p)(I_1-I)K\sum_{j=0}^\infty\pmatrix{-1/2
\cr j}
(-1)^jL^j.\nonumber
\end{eqnarray}

The following lemma is used as the main tool in the proofs.
%%%%%%%%%%%%%%%%%%%%%%%%%%%%%%%%%%%%%%%%%%%%%%%%%%%%%%%%%%%%%%
%
%l
\begin{lemma}\label{main}
If condition~\textup{(\ref{cond1})} is satisfied, then
%e41 ###
%e40 ###
%e39 ###
%e38 ###
%e37 ###
%e36 ###
%
\begin{eqnarray}
\Lambda_1&=& I+\sum_{j=1}^3a_jY^j+C\overline{q}^3(p+\overline
{q})Y^4\Theta,\label{L1}\\
\ln{\Lambda_1}&=&\sum_{j=1}^3\gamma_jY^j+C\overline
{q}^3(p+\overline{q})Y^4\Theta, \label{lnL1}\\
\ln{\Lambda_1}&=&\gamma_1Y+\frac{19}{60}\gamma_1Y^2\Theta,\label
{shlnL1}\\
\|\Lambda_2\|&\leq&\frac{19}{30},\qquad \|\Lambda_1-I\|\leq0.1,\label
{L12}\\
W_1&=&I+(\varkappa_1+\varkappa_2)Y+C\overline{q}(p+\overline
{q})Y^2\Theta,\label{W1}\\
W_2&=&C(p+\overline{q})(I_1-I)\Theta,\qquad \|W_2\|\leq7. \label{W2}
\end{eqnarray}
For any finite signed measure $M$ on $\mathbb{Z}$ and any $t>0$, we have
%e43 ###
%e42 ###
%
\begin{eqnarray}
\|M\exp\{t\ln\Lambda_1\}\|&\leq&
C\|M\exp\{(t\gamma_1/30)Y\}\|,\label{smL1}
\\
\|MD_j^t\|&\leq& C\|M\exp\{(t\gamma_1/30)Y\}\|, \qquad j=1,2,3.\label{smD}
\end{eqnarray}
Estimates (\ref{smL1})--(\ref{smD}) also hold for the local norm.
% if the total variation norm on both sides is replaced by the
%local one.
\end{lemma}

%%%%%%%%%%%%%%%%%%%%%%%%%%%%%%%%%%%%%%%%%%%%%%%%%%%%%%%%%%%%%%%%%%%%%%%
%
\begin{pf}
We have
\[
a_1=\gamma_1\leq\frac{1}{30},\qquad
\frac{1}{q+\overline{q}}\leq\frac{1}{1-p}\leq2,\qquad
\|Y\|\leq\|G\| +1=2,
\]
%
%e44 ###
%
\begin{eqnarray}
\|K\| &\leq&\sum_{j=0}^\infty(p+\overline q+4 a_1)^j\leq 3,\label{normK}
\\
\|L\| &\leq& 9\cdot4\cdot\overline q^2\cdot 4 \biggl(1+\frac{p}{q+\overline{q}}(1+p) \biggr)\leq 0.4.\label{normL}
\end{eqnarray}
Note that
\[
\left\vert\pmatrix{1/2 \cr2} \right\vert=\frac{1}{8},\qquad
\left\vert\pmatrix{1/2 \cr3} \right\vert=\frac{1}{16},\qquad
\left\vert\pmatrix{1/2 \cr j} \right\vert\leq\frac{5}{128},\qquad
j\geq4.
\]
We have
\[
\sum_{j=1}^{\infty} \left\vert\pmatrix{1/2 \cr j} \right\vert\|L\|^{j-1}\leq
\frac{1}{2}+\frac{0.4}{8}+\frac{(0.4)^2}{16}+\frac{5}{128}\frac{(0.4)^3}{0.6}\leq 0.5642
\]
and
\begin{eqnarray*}
[I(1+\overline q)-pI_1+2a_1Y]L
&=& 4\frac{\overline q^2 Y^2}{(q+\overline{q})^2}[q^2I+p(q+\overline{q})(I-pI_1)]K\nonumber\\
&=&\gamma_1Y^2\frac{12\overline q}{q(q+\overline{q})}\bigl(q^2+p(q+\overline{q})(1+p)\bigr)\Theta=\gamma_1Y^2\Theta.\nonumber
\end{eqnarray*}
Consequently,
\[
\Lambda_1=I+\gamma_1Y+\tfrac{1}{2}0.5642\gamma_1Y^2\Theta=I+
1.2821\gamma_1Y\Theta=
I+0.1\Theta
\]
and
\begin{eqnarray*}
\ln\Lambda_1
&=&\Lambda_1-I+\sum_{j=2}^{\infty}\frac{(-1)^{j+1}}{j}(\Lambda_1-I)^j\nonumber\\
&=&\gamma_1Y+0.2821\gamma_1Y^2\Theta+\frac{1}{2}1.2821\gamma_1^2Y^2\sum_{j=2}^{\infty}(0.1)^{j-2}\Theta\ =\gamma_1Y+\frac{19}{60}\gamma_1Y^2\Theta.
\end{eqnarray*}
Moreover, $\|\Lambda_2\|\leq
p+\overline{q}+\|\Lambda_1-I\|\leq19/30$.
Thus, we have proven
(\ref{shlnL1}). We use this estimate for obtaining (\ref{smL1}).
By the properties of the total variation norm, we have
\[
\|M \mathrm{e}^{t \ln \Lambda_1}\|
\leq\bigg\|M \exp\biggl\{\frac{t\gamma_1}{30}Y \biggr\} \bigg\|
\bigg\|\exp\biggl\{\frac{29t\gamma_1}{30}Y+\frac{19t\gamma_1}{60}Y^2\Theta\biggr\} \bigg\|.
\]
Applying Lemma~\ref{smooth}, we prove that the second norm is majorized by
\[
1+\sum_{r=1}^{\infty}\frac{1}{r!}\bigg\|\frac{19}{60}t\gamma_1 Y^2 \exp\biggl\{\frac{29}{30r}t\gamma_1Y \biggr\} \bigg\|^r\leq1+
\sum_{r=1}^{\infty}\frac{\mathrm{e}^r}{r^r \sqrt{2\curpi
r}} \biggl(\frac{57r}{58
\mathrm{e}} \biggr)^r\leq C.
\]
\\
The last two estimates obviously lead to (\ref{smL1}). The estimate
(\ref{smD}) is proved similarly. For the proof of (\ref{L1}), note
that
%e50 ###
%e49 ###
%e48 ###
%e47 ###
%e46 ###
%
\begin{eqnarray}
\Lambda_1&=&I+\gamma_1Y-\frac{\overline
q^2}{(q+\overline{q})^2}Y^2[q^2I+ p(q+\overline{q})(I-pI_1)]K+C\overline
q^4 Y^4
\Theta,\label{L01}\\
I_1-I&=&Y(I-pI_1),\qquad
(q+\overline{q})B=I+p(I_1-I)B,\\
I_1-I&=&2Y\Theta,\qquad (q+\overline{q})(I-pI_1)B=qI-p\overline q
(I-pI_1)BY,\label{L02}\\
B&=&\frac{1}{q+\overline{q}}I+\frac{pq}{(q+\overline{q})^2}Y-\frac
{p^2\overline
q}{(q+\overline{q})^2}Y(I_1-I)B\label{L4},\\
K&=&B+2\gamma_1Y
KB=B-2\gamma_1B^2Y+4\gamma_1^2B^2Y^2B^2K.\label{L3}
\end{eqnarray}
Substituting
(\ref{L02})--(\ref{L3}) into (\ref{L01}), we obtain (\ref{L1}).
Taking into account (\ref{L12}), we obtain
\[
\ln\Lambda_1=\sum_{j=1}^{3}\frac{(-1)^{j+1}}{j}(\Lambda
_1-I)^j+C(\Lambda_1-I)^4\Theta.
\]
Now, for the proof of (\ref{lnL1}), it suffices to use (\ref{L1}).
From (\ref{L3}) and the first relation in (\ref{L02}), we get
$p(I_1-I)K=p(I_1-I)B+Cp\overline q Y^2\Theta$.
Moreover,
\[
p(I_1-I)B=\frac{pq}{q+\overline{q}}Y-\frac{p^2\overline
q(I_1-I)BY}{q+\overline{q}}=\frac{pqY}{q+\overline{q}}+Cp\overline q
Y^2\Theta.
\]
The last two equations and (\ref{a13}) allow us to
prove (\ref{W1}). Since $W_1+W_2=I$, we easily obtain the
first relation in (\ref{W2}). Now,
\[
\|W\|_2\leq\frac{1}{2} \Biggl(1+(q+\overline{q})3\sum_{j=0}^{\infty}
\left\vert\pmatrix{-1/2 \cr j} \right\vert
0.4^j \Biggr) +|\overline{q}-p|\cdot2\cdot
3\sum_{j=0}^{\infty} \left\vert\pmatrix{-1/2 \cr j} \right\vert0.4^j<7.
\]
Thus, Lemma~\ref{main} is proved. For the lower bound estimates, we
need the following result.
\end{pf}

%l5.5
\begin{lemma}\label{lbinvlemma}
Let M be concentrated on $\mathbb{Z}$, $\alpha\in\mathbb{R}$ and
$b>1$. Then,
%e52 ###
%e51 ###
%
\begin{eqnarray}\label{lbinv}
\|M\|&\geq& C
\bigg\vert\int_{-\infty}^{\infty}\mathrm{e}^{-{t^2}/{2}}\widehat
M \biggl(\frac{t}{b} \biggr)\mathrm{e}^{-\mathrm{i}t\alpha}\,{\mathrm{d}}t \bigg\vert
, \label{lb1}\\
\|M\|_\infty&\geq& \frac{C}{b} \bigg\vert
\int_{-\infty}^{\infty}\mathrm{e}^{-{t^2}/{2}}\widehat
M \biggl(\frac{t}{b} \biggr)\mathrm{e}^{-\mathrm{i}t\alpha}\,{\mathrm{d}}t \bigg\vert.
\label{lb2}
\end{eqnarray}
The estimates (\ref{lb1}) and (\ref{lb2}) remain valid if
$\mathrm{e}^{-{t^2}/{2}}$ is replaced by $t \mathrm{e}^{{-t^2}/{2}}.$
\end{lemma}

Lemma~\ref{lbinvlemma}, with $\|M\|$ replaced by the uniform
norm of $M$, was proven in \cite{p4}. Since the uniform norm is
majorized by the total variation norm, (\ref{lb1}) also holds.
%
%l5.6
\begin{lemma}\label{charf}
If (\ref{cond1}) is satisfied, then, for all $|t|\leq\curpi$,
%e54 ###
%e53 ###
%
\begin{eqnarray}
&\big|\exp\bigl\{n\varkappa_1\bigl(\widehat Y(t)-\mathrm{i}t/q\bigr)\bigr\}-1\big|\leq C n\overline{q}t^2,&\label{ch1}\\
&|\widehat{D}_2^n(t)|\leq1,\qquad |\widehat{D}_2^n \exp\{-\mathrm{i}tn\gamma_1/q\}-1| \leq Cn\overline{q}t^2.& \label{ch2}
\end{eqnarray}
\end{lemma}

Proof of Lemma~\ref{charf} is straightforward and therefore
omitted.

Finally, let us introduce an inverse compound measure for $H$. Let
\[
H^{-1}=\exp\Biggl\{-\sum_{j=1}^{\infty}\frac{p^j}{j} \biggl(1- \biggl(\frac{1-p_0q/(q+\overline{q})}{1-\varkappa_2} \biggr)^j \biggr) (I_j-I)\Biggr\}.
\]
%
%l5.7
\begin{lemma}
If (\ref{cond1}) is satisfied, then
%e55 ###
%
\begin{eqnarray}\label{h}
\|H^{-1}\|\leq\mathrm{e}^2
\end{eqnarray}
and, for any (signed) finite measure $M$ concentrated at $\mathbb{Z}$,
%e56 ###
%
\begin{equation}\label{hm}
\|MH\|\geq\mathrm{e}^{-2}\|M\|,\qquad
\|M\exp\{-p_0\gamma_1Y\}\|\geq\|M\|.
\end{equation}
The estimates in (\ref{hm}) remain valid if the total variation
norm is replaced by the local norm.
\end{lemma}

\begin{pf}
Estimate (\ref{h}) easily follows from the property $\|\mathrm{e}\|^M\leq\mathrm{e}^{\|M\|}$; see the \hyperref[sec1]{Introduction}.
Now,
$\|M\|=\|MHH^{-1}\|\leq\|MH\| \|H^{-1}\|\leq\mathrm{e}^2\|MH\|$.
Since $\exp\{p_0\gamma_1Y\}$ is a distribution, its total
variation is 1. Therefore,
$\|M\|=\|M\exp\{-p_0\gamma_1Y\}\exp\{p_0\gamma_1Y\}\|\leq
\|M\exp\{-p_0\gamma_1Y\}\|$. Estimates for the local norm are proved
similarly.
\end{pf}
%
%ADD HERE
%%%%%%%%%%%%%%%%%%%%%%%%%%%%%%%%%%%%%%%%%%%%%%%%%%%%%%%%%%%%%%%%%%%%%%
%%%%%%%%%%%%%%%%%%%%%%%%%%%%%%%%%%%%%%%%%%%%%%%%%%%%%%%%%%%%%%%%%%%%%%
%s6 ###
\section{Proofs}
For upper bound estimates, we use an adaptation of Le~Cam's
\cite{Lec60}
approach which deals with convolutions of measures.
\begin{pf*}{Proof of Theorem~\ref{T1}} Without loss of generality, we
can assume
that (\ref{cond1}) holds. We have
\[
\|F_n-HD_1^\lambda\|\leq\|\Lambda_1^n-D_1^n\|\|W_1\|+\|D_1^n(W_1-H\exp\{-p_0\gamma_1Y\})\|+\|\Lambda_2\|^n\|W_2\|.
\label{p1}
\]
Further, in view of Lemma~\ref{main},
\begin{eqnarray*}
\|\Lambda_1^n-D_1^n\|&\leq&
\bigg\|D_1^n\int_0^1(\exp\{\tau[n\ln\Lambda_1-n\gamma_1Y]\})_\tau'\,{\mathrm{d}}\tau \bigg\|\\
&\leq& n\int_0^1\|[\ln\Lambda_1-\gamma_1Y]\exp\{\tau n\ln
\Lambda_1+(1-\tau)n\gamma_1Y\}\|\,{\mathrm{d}}\tau\\
&\leq&
Cn\|[\ln\Lambda_1-\gamma_1Y]\exp\{(n\gamma_1/30)Y\}\|\leq
Cn\overline{q}^2\|Y^2\exp\{(n\gamma_1/30)Y\}\|.
\end{eqnarray*}
By Lemma~\ref{main},
\begin{eqnarray*}
W_1-H\exp\{-p_0\gamma_1Y\}&=&[W_1-I-(\varkappa_1+\varkappa_2)Y]+
[I+(\varkappa_1+\varkappa_2)Y-(I-p_0\gamma_1Y)H]\nonumber\\
&&{}+[H(I-p_0\gamma_1Y-\exp\{-p_0\gamma_1Y\})]=C\overline
{q}(p+\overline{q})Y\Theta.
\label{pp1}
\end{eqnarray*}
Taking into account the last two estimates, applying Lemma
\ref{smooth} and estimating $\|W_2\|$ and $\|\Lambda_2\|$ by
(\ref{W2}) and (\ref{L12}), we complete the proof of (\ref{CP1}).
The estimates in (\ref{CP2}) and (\ref{CP3}) are proved similarly.
\end{pf*}
\begin{pf*}{Proof of Corollary ~\ref{cor2}} Following the proof of
(\ref{smL1}), one can prove the same property for $H_1$. Also,
\[
\|HD_1^\lambda-HH_1^\lambda\| \leq
C\lambda\|(D_1-H_1)\exp\{(n\gamma_1/30)Y\}\|\leq
Cn\overline{q}^2\|Y^2\exp\{(n\gamma_1/30)Y\}\|.
\]
The rest of the proof is obvious.
\end{pf*}
\begin{pf*}{Proof of Theorem~\ref{T2}} We have
\begin{eqnarray*}
&&\|F_n-H D_1^\lambda(I+n\gamma_2Y^2)\|
\\
&&\quad\leq
\|\Lambda_2\|^n\|W_2\|+\|W_1\| \|\Lambda_1^n-D_2^n\|\\
&&\qquad{}+\|W_1\| \|D_1^n(\mathrm{e}^{n\gamma_2Y^2}-I-n\gamma_2Y^2)\|+
\|D_1^n(I+n\gamma_2Y^2)(W_1-H\mathrm{e}^{-p_0\gamma_1Y})\|.
\end{eqnarray*}
%
%%%%%%%%%%%%%%%%%%%%%%%%%%%%%%%%%%%%%%%%%%%%%%%%%%%%%%%%%%%%%%%%%%%%%%
Similarly to the proof of Theorem~\ref{T1}, and using (\ref{smD}), we
obtain
\begin{eqnarray*}
\|\Lambda_1^n-D_2^n\|&\leq&
Cn\bigg\|[\ln\Lambda_1-\gamma_1Y-\gamma_2Y^2]\int_0^1 \exp\{\tau n\ln
\Lambda_1+(1-\tau)[n\gamma_1Y+n\gamma_2Y^2]\}\,{\mathrm{d}}\tau \bigg\|\\
&\leq&Cn\|[\ln\Lambda_1-\gamma_1Y-\gamma_2Y^2]\exp\{(n\gamma
_1/30)Y\}\|\\
&\leq&
Cn\overline{q}^2(\overline{q}+p)\|Y^3\exp\{(n\gamma_1/30)Y\}\|
\end{eqnarray*}
and
\begin{eqnarray*}
\|D_1^n(\mathrm{e}^{n\gamma_2Y^2}-I-n\gamma_2Y^2) \|
&\leq&\bigg\|(n\gamma_2Y^2)^2\int_0^1D_1^n\mathrm{e}^{\tau n\gamma_2Y^2}(1-\tau)\,{\mathrm{d}}\tau \bigg\|
\\
&\leq& C(n\gamma_2)^2\|Y^4\exp\{(n\gamma_1/30)Y\}\|.
\end{eqnarray*}
Note that for any signed finite measure $M$,
\begin{eqnarray*}
\|D_1^n(I+n\gamma_2Y^2)M\|
\leq \|D_1^{n/2}\|M(1+n|\gamma_2|\|Y^2 D_1^{n/2}\|)
\leq C\|D_1^{n/2}M\|.
\end{eqnarray*}
The rest of the proof is very similar to the proof of Theorem
\ref{T1} and is hence omitted.
\end{pf*}
\begin{pf*}{Proof of Theorems~\ref{T3},~\ref{T4} and~\ref{T5}} The
proofs are very similar to those of Theorems~\ref{T1} and~\ref{T2}.
From Lemma
\ref{main} and the definition of the exponent measure, it is not
difficult to show that
\begin{eqnarray*}
W_1-\mathrm{e}^{\varkappa_1Y}H&=&[W_1-I-(\varkappa_1+\varkappa_2)Y]+
[I+(\varkappa_1+\varkappa_2)Y-(I+\varkappa_1Y)H]\\
&&{}+H(I+\varkappa_1Y-\mathrm{e}^{\varkappa_1Y})=C\overline
{q}(p+\overline{q})Y^2\Theta,\\
\|\Lambda_1^nW_1-D_2^nH\mathrm{e}^{\varkappa_1Y}\|&\leq&
\|\Lambda_1^n-D_2^n\| \|W_1\|+\|D_2^n(W_1-H\mathrm{e}^{\varkappa
_1Y})\|.
\end{eqnarray*}
Now, it is not difficult to prove Theorem~\ref{T3}. Theorem~\ref{T5}
is proved similarly. For the proof of Theorem~\ref{T4}, one should
use Theorem~\ref{T5}, the triangle inequality and the fact that
\begin{eqnarray*}
\|D_2^n(I+n\gamma_3Y^3)-D_3^n\|&=&\bigg\|D_2^n\int_0^1(1-\tau)\mathrm
{e}^{\tau n\gamma_3Y^3}(n\gamma_3Y^3)^2\,{\mathrm{d}}\tau \bigg\|\\
&\leq&C(n\gamma_3)^2\|Y^6\exp\{n\gamma_1Y^2/30\}\|.
\end{eqnarray*}
For the last estimate, we have used the same argument as in the
proof of (\ref{smD}).
\end{pf*}
\begin{pf*}{Proof of Theorem~\ref{T6}} Taking into account Theorem
\ref{T2}, (\ref{a67}) and (\ref{hm}), we get
%e57 ###
%
\begin{eqnarray}\label{pl1}
\|F_n-HD_1^\lambda\|&\geq&
n|\gamma_2|\|HD_1^\lambda Y^2\|-C\overline{q}\biggl(\overline{q}+\frac
{p}{\sqrt{n\overline{q}}} \biggr)\nonumber\\[-8pt]\\[-8pt]
&\geq&
C_{10}n|\gamma_2|\|D_1^n(I_1-I)^2\|-C_{11}\overline{q}\biggl(\overline
{q}+\frac{p}{\sqrt{n\overline{q}}} \biggr).\nonumber
\end{eqnarray}
Let $z=t/(h\sqrt{n\overline{q}})$ and $\mu=n\gamma_1/q$. The constant
$h>1$ will be chosen later. Applying Lemma~\ref{charf}, we
then obtain
%e58 ###
%
\begin{equation}
J=\bigg\vert\int_\mathbb{R}\mathrm{e}^{-t^2/2}\widehat D_1^n(z)\mathrm{e}^{-\mathrm{i} z\mu}(\mathrm{e}^{\mathrm{i}z}-1)^2\,{\mathrm{d}} t
\bigg\vert\geq\bigg\vert\int_\mathbb{R}\mathrm{e}^{-t^2/2}z^2\,{\mathrm{d}}t\bigg\vert-J_1-J_2.
\label{pl2}
\end{equation}
Here,
\begin{eqnarray*}
J_1&=&\int_\mathbb{R}\mathrm{e}^{-t^2/2}z^2|\widehat D_1^n(z)\mathrm
{e}^{-\mathrm{i} z\mu}-1|\,{\mathrm{d}}t\leq Cn\overline{q}\int
_\mathbb{R}z^4\mathrm{e}^{-t^2/2}\,{\mathrm{d}}t=\frac{C}{h^4n\overline{q}},
\\
J_2&=&\int_\mathbb{R}\mathrm{e}^{-t^2/2}|\widehat D_1^n(z)\mathrm
{e}^{-\mathrm{i} z\mu}||(\mathrm{e}^{\mathrm{i}z}-1)^2-(\mathrm
{i}z)^2|\,{\mathrm{d}}t\leq
\frac{C}{h^3n\overline{q}\sqrt{n\overline{q}}}.
\end{eqnarray*}
Combining the last two
estimates with (\ref{pl2}) and choosing $h$ to be a sufficiently
large absolute constant, we obtain
\[
J\geq\frac{C_{12}}{h^2n\overline{q}}
\biggl(1-\frac{C_{13}}{h^2}-\frac{C_{14}}{h\sqrt{n\overline{q}}} \biggr)\geq
\frac{C_{15}}{n\overline{q}}.
\label{pl3}
\]
Applying Lemma~\ref{lbinvlemma} and substituting the result into
(\ref{pl1}), we get (\ref{LB1}). Estimates (\ref{LB2}) and
(\ref{LB3}) are proved similarly.
%%%%%%%%%%%%%%%%%%%%%%%%%%%%%%%%%%%%%%%%%%%%%%%%%%%%%%%%%%%%%%%%%%%%%%%%%%%%%

For the proof of Theorem~\ref{T7}, one should use Theorem~\ref{T4}
and take $t\exp\{-t^2/2\}$ instead of $\exp\{-t^2/2\}$. The proof is
then almost identical to that of Theorem~\ref{T6} and is
hence omitted.
\end{pf*}
\begin{pf*}{Proof of Theorem~\ref{T8}} We have
\begin{eqnarray*}
\big\vert\|F_n-HD_1^\lambda\|-A_{11} \big\vert
&\leq&\|F_n-HD_1^\lambda(I+n\gamma_2Y^2)\|
\\
&&{}+\|(H\mathrm{e}^{-p_0\gamma_1Y}-I)D_1^nn\gamma_2Y^2 \|
+n|\gamma_2|\bigg\|\biggl(Y^2-\frac{1}{q^2}(I_1-I)^2 \biggr)D_1^n \bigg\|\\
&&{}+\frac{n|\gamma_2|}{q^2}\bigg\|(I_1-I)^2 \biggl(D_1^n-\exp\biggl\{\frac{n\gamma
_1}{q}(I_1-I) \biggr\} \biggr) \bigg\|\\
&&{}+
\bigg\vert\frac{n|\gamma_2|}{q^2}\bigg\|(I_1-I)^2\exp\biggl\{\frac{n\gamma_1}{q}(I_1-I) \biggr\} \bigg\|-A_{11} \bigg\vert.
\end{eqnarray*}
One should now apply Theorem~\ref{T2}, (\ref{L02}), Lemmas
\ref{sharpC},~\ref{smooth} and the following, easily verifiable, relations:
\[
Y=\frac{(I_1-I)}{q}\sum_{j=0}^\infty\biggl(\frac{p}{q} \biggr)^j(I_1-I)^j=
\frac{(I_1-I)}{q}+\frac{3p}{q^2}(I_1-I)^2\Theta
\]
and
\begin{eqnarray*}
&&D_1-\exp\biggl\{\frac{\gamma_1}{q}(I_1-I) \biggr\}\\
&&\quad=\exp\biggl\{\frac{\gamma_1}{q}(I_1-I) \biggr\} \biggl(
\exp\biggl\{\frac{3p\gamma_1}{q^2}(I_1-I)^2\Theta \biggr\}-I\biggr)=
Cp\overline{q}(I_1-I)^2\Theta.
\end{eqnarray*}
Note that
\begin{eqnarray*}
&&\bigg\|(I_1-I)^2 \biggl(D_1^n-\exp\biggl\{\frac{n\gamma_1}{q}(I_1-I) \biggr\} \biggr) \bigg\|
\\
&&\quad=
\Bigg\|(I_1-I)^2 \biggl(D_1-\exp\biggl\{\frac{\gamma_1}{q}(I_1-I) \biggr\} \biggr) \sum
_{j=1}^nD_1^{n-j}\exp\biggl\{(j-1)\frac{\gamma_1}{q}(I_1-I) \biggr\} \Bigg\|\\
&&\quad\leq Cnp\overline{q}\biggl(
\|(I_1-I)^4D_1^{n/3} \|+\bigg\|(I_1-I)^4\exp\biggl\{\frac{n\gamma_1}{3q}(I_1-I)\biggr\} \bigg\|\biggr)
.
\end{eqnarray*}
All other estimates are obtained similarly.
\end{pf*}

\section*{Acknowledgements}

The main part of this work was accomplished during the first
author's stay at the Department of Mathematics, IIT Bombay, during
November 2006. The first author would like to thank the members
of the Department for their hospitality. Also, we are grateful to the
referees for some useful remarks.

\printhistory

\end{document}